\documentclass[12pt]{article}
\usepackage{amsmath,amsfonts,amsthm,amscd}
\newtheorem{theorem}{Theorem}

\newtheorem{proposition}{Proposition}
\newtheorem{lemma}{Lemma}
\newtheorem{definition}{Definition}
\newtheorem{corollary}{Corollary}
\theoremstyle{remark}
\newtheorem*{remark}{Remark}

\numberwithin{equation}{section}
\numberwithin{theorem}{section}
\numberwithin{proposition}{section}
\numberwithin{lemma}{section}
\numberwithin{claim}{section}
\numberwithin{corollary}{section}
\begin{document}
\title{A fully nonlinear equation on four-manifolds 
with positive scalar curvature}
\author{ {Matthew J. Gursky}\thanks{Supported in part by NSF 
Grants DMS-0200646 and INT-0229457.} \and 
{Jeff A. Viaclovsky}\thanks{Supported in
part by NSF Grant DMS-0202477.}}
\date{January 30, 2003}
\maketitle
\begin{abstract}
We present a conformal deformation involving a fully nonlinear
equation in dimension $4$, starting with a metric of positive 
scalar curvature. Assuming a certain conformal 
invariant is positive, one may deform from 
positive scalar curvature to a stronger 
condition involving the Ricci tensor.
A special case of this deformation provides an alternative 
proof to the main result in \cite{CGY1}. We also give 
a new conformally invariant condition for positivity of the 
Paneitz operator, generalizing the results in \cite{Gursky1}.
From the existence results in \cite{ChangYang}, this allows us
to give many new examples of manifolds admitting metrics with 
constant $Q$-curvature.
\end{abstract}
\section{Introduction}

Let $(M,g)$ denote a closed, 4-dimensional Riemannian 
manifold, and let $Y[g]$ denote the Yamabe 
invariant of the conformal class $[g]$: 
\begin{align}
Y[g] \equiv \underset{ \tilde{g} \in [g] }{\mbox{inf }} Vol(\tilde{g})^{-1/2} \int_M R_{\tilde{g}}
 dvol_{\tilde{g}},
\end{align}
where $R_{\tilde{g}}$ denotes the scalar curvature.
Another important conformal invariant is 
\begin{align}
\label{F2}
\mathcal{F}_2([g])\equiv
\int_M \left( -\frac{1}{2} |Ric_g|^2 + \frac{1}{6} R_g^2 \right) dvol_g,
\end{align}  
where $Ric_g$ is the Ricci tensor. By the Chern-Gauss-Bonnet formula (\cite{Besse}),
\begin{align}
\label{CGB}
8\pi^2\chi(M) = \int_M |W_g|^2 dvol_g  + \mathcal{F}_2([g]).
\end{align}
Thus, the conformal invariance of $\mathcal{F}_2$ follows from the well known
(pointwise) conformal invariance of the Weyl tensor $W_g$ 
(see \cite{Eisenhart}).  

Define the tensor
\begin{align}
A^t_g = \frac{1}{2} \Big( Ric_g - \frac{t}{6} R_g g \Big).
\end{align}
Note that for $t=1$, 
$A^1_g$ is the classical Schouten tensor (\cite{Eisenhart}).
Let $\sigma_2(g^{-1}A^t_g)$ 
denote the second elementary symmetric 
function of the eigenvalues of $g^{-1}A^t_g$, viewed as an endomorphism of the tangent bundle. 
Then a simple calculation gives
\begin{align}
\mathcal{F}_2([g]) = 4 \int_M \sigma_2(g^{-1}A^1_g) dvol_g. 
\end{align}  
Our main result is the following:
\begin{theorem}
\label{main}
Let $(M,g)$ be a closed 4-dimensional Riemannian 
manifold with positive scalar curvature. If  
\begin{align} 
\label{assump}
 \mathcal{F}_2([g]) + \frac{1}{6} (1-t_0) (2-t_0) (Y[g])^2 > 0,
\end{align}
for some $t_0 \leq 1$, 
then there exists a conformal metric $\tilde{g} = e^{-2u} g$ with
$R_{\tilde{g}} > 0$ and $\sigma_2(A^{t_0}_{\tilde{g}}) > 0$ pointwise. 
This implies 
the pointwise inequalities
\begin{align}
\label{Ricineq}
(t_0-1) R_{\tilde{g}} \tilde{g} < 2 Ric_{\tilde{g}} < (2-t_0) 
R_{\tilde{g}} \tilde{g}.
\end{align}
\end{theorem}

As applications of Theorem \ref{main}, we consider two different 
values of $t_0$.  When $t_0 = 1$, we obtain 
a different proof of the following result in \cite{CGY1}:
\begin{corollary}
\label{AMP}
Let $(M,g)$ be a closed 4-dimensional Riemannian 
manifold with positive scalar curvature. If
$\mathcal{F}_2([g]) > 0$, 
then there exists a conformal metric $\tilde{g} = e^{-2u} g$ with
$R_{\tilde{g}} > 0$ and $\sigma_2(\tilde{g}^{-1}A^1_{\tilde{g}}) > 0$ pointwise.  In particular, the Ricci curvature of 
$\tilde{g}$ satisfies 
$$
0 < 2 Ric_{\tilde{g}} <  R_{\tilde{g}} \tilde{g}.
$$
\end{corollary}

The proof in \cite{CGY1} involved regularization by a fourth-order 
equation and relied on some delicate integral estimates.  
By contrast, the proof of Theorem \ref{main} seems more direct, and 
depends on general {\it a priori} estimates for fully nonlinear equations  
developed in \cite{Jeff1}, \cite{GuanWang1}, \cite{LiLi2}, 
and \cite{GVNegative}.

Our second application is to the spectral properties of a conformally invariant differential operator known as the 
{\it Paneitz operator}.  Let $\delta$ denote the $L^2$-adjoint of the exterior derivative $d$; then the Paneitz operator is
defined by 
\begin{align}
P_{g} \phi   = \Delta^2 \phi + \delta \Big( \frac{2}{3} R_g g - 2 Ric_g \Big)
d\phi. 
\end{align}
The Paneitz operator is conformally invariant, in the sense
that if $\tilde{g} = e^{-2u}g$, then
\begin{align}
\label{confinv}
P_{\tilde{g}} = e^{4u}P_g.
\end{align}
Since the volume form of the conformal metric $\tilde{g}$ is $dvol_{\tilde{g}} = e^{-4u}dvol_g$, an immediate consequence of 
(\ref{confinv}) is the conformal invariance of the Dirichlet energy 
\begin{align*}
\langle P_{\tilde{g}} \phi, \phi \rangle_{L^2(M,\tilde{g})} = \langle P_g \phi, \phi \rangle_{L^2(M,g)}.
\end{align*}
In particular, positivity of the Paneitz operator is a conformally invariant property, and clearly the kernel is invariant as well.

To appreciate the geometric significance of the Paneitz operator we need to define the 
associated $Q$-{\it curvature}, introduced by Branson: 
\begin{align} 
Q_g = - \frac{1}{12} \Delta R_g  + 2 \sigma_2(g^{-1}A^1_g).
\end{align}
Under a conformal change of metric 
$\tilde{g} = e^{-2u}g$, the $Q$-curvature transforms according to the 
equation
\begin{align}
\label{Qchange}
-Pu + 2Q_g = 2Q_{\tilde{g}}e^{-4u},
\end{align}
see, for example, \cite{BransonOrsted}.
Note that 
\begin{align}
\label{intQ}
\int_M Q_g dvol_g = \frac{1}{2} \mathcal{F}_2([g]),
\end{align}
so the integral of the $Q$-curvature is conformally invariant.

The $Q$-curvature and Paneitz operator have become important objects of study in the geometry of four-manifolds, 
and play a role in such diverse topics as Moser-Trudinger inequalities (\cite{Beckner}, \cite{BransonChangYang}),
compactification of complete conformally flat manifolds (\cite{ChangQingYang}), twistor theory (\cite{ES}),
gauge choices for Maxwell's equations (\cite{EastwoodSinger2}),
and most recently in the study of conformally compact $AHE$ manifolds (\cite{FeffermanGraham}, \cite{GrahamZworski}).

Our interest here is in the spectral properties of the Paneitz operator and the related question of the existence of metrics with
constant $Q$-curvature.  The most general work on this subject was done by Chang and Yang
\cite{ChangYang}, who studied the problem of constructing 
conformal metrics with constant $Q$-curvature by
variational methods.  They considered the functional 
\begin{align}
\label{FP}
F[\phi] = \langle P_g \phi, \phi \rangle - 4\int_M Q\phi dvol - \left( \int_M Q dvol \right) \log \int_M e^{-4\phi} dvol,
\end{align}
and analyzed the behavior of a minimizing sequence.  Of course, it is not clear {\it a priori} that 
$F$ is even bounded from below.  Indeed, if the Paneitz operator has a negative eigenvalue and the conformal invariant
(\ref{intQ}) is positive, then Chang and Yang showed that $\inf F = -\infty$ (see \cite{ChangYang}, p. 177).  For example,  
take a compact surface $\Sigma$ of curvature $-1$ with first eigenvalue $\lambda_1(-\Delta) << 1$.  Then the product manifold
$M = \Sigma \times \Sigma$ will have $\lambda_1(P) < 0$ and $\int Q dvol > 0$.  

Chang and Yang also pointed out the connection between 
the conformal invariant (\ref{intQ}) and the best constant in
the inequality of Adams \cite{Adams}, another key point for establishing the 
$W^{2,2}$ compactness of a minimizing sequence.  
Combining these observations, they proved:

\begin{theorem} (\cite{ChangYang})  Let $(M,g)$ be a compact 4-manifold.
Assume $(i)$ the Paneitz 
operator $P_g$ is non-negative with $Ker P = \{constants\}$, and $(ii)$ the conformal
invariant (\ref{intQ}) is strictly less than the value attained by
 the round sphere.  Then there exists a minimizer
of $F$, which satisfies (\ref{Qchange}) with $Q_{\tilde{g}} = constant$.
\end{theorem}
Subsequently, the first author proved that any four-manifold of positive scalar curvature which is not conformally equivalent to the sphere already satisfies the
second assumption of Chang-Yang.  In addition,

\begin{theorem} (\cite{Gursky1}) 
Let $(M,g)$ be a compact 4-manifold.  If the scalar curvature of $g$ is non-negative
and $\int Q dvol \geq 0$, then the Paneitz operator is positive and $Ker P = \{constants\}$.
\end{theorem}

Because of the example of Chang-Yang, it is clear that one cannot relax the condition on the scalar curvature in the
above theorem.  On the other hand, the positivity of the conformal invariant (\ref{intQ}) is a rather strong assumption.
For example, if the scalar curvature is strictly positive, then the positivity of (\ref{intQ}) implies the 
vanishing of the first Betti number of $M$ (see \cite{Gurskyweyl}).  Thus, for example, the manifold $N \# (S^1 \times S^3)$ 
can not admit a metric of positive 
scalar curvature with $\int Q dvol > 0$.

It is interesting to note that the positivity of the Paneitz operator
 was studied by Eastwood and Singer in \cite{ES} for reasons motivated
by twistor theory.  They constructed metrics on $k(S^3 \times S^1)$ 
for all $k>0$ with $P \geq 0$ and $Ker P = \{constants\}$.  Since these manifolds have $\int Q dvol < 0$, the Eastwood-Singer construction is in some respects complementary to the 
result of \cite{Gursky1}.

By combining Theorem \ref{main} with
$t_0=0$, and an integration by parts argument, 
we obtain a new criterion for the positivity of $P$:

\begin{theorem}
\label{Q}

Let $(M,g)$ be a closed 4-dimensional Riemannian 
manifold with positive scalar curvature. If
\begin{align} 
\label{assump3}
 \int Q_g dvol_g + \frac{1}{6}(Y[g])^2 > 0,
\end{align}
then the Paneitz operator is nonnegative, and 
$Ker P = \{ constants \}$. Therefore, 
by the results in \cite{ChangYang}, there 
exists a conformal metric $\tilde{g} = e^{-2u} g$
with $Q_{\tilde{g}} = constant.$
\end{theorem}

Since Theorem \ref{Q} allows the integral of the $Q$-curvature to be negative, we are able to use surgery techniques to construct many new examples of 
manifolds which admit metrics with constant $Q$. 
For example, we will show that 
\begin{align*}
N &=  (S^2 \times S^2) \# k(S^1 \times S^3), \ k \leq 5,\\
N &=  \mathbb{CP}^2 \# k(S^1 \times S^3), \ k \leq 5, \\ 
N &= \mathbb{CP}^2 \# k (\mathbb{RP}^4), \ k \leq 8, \\ 
N &= k(S^1 \times S^3) \#  l(\mathbb{RP}^4), \ 2k +l \leq 9,
\end{align*} 
all admit metrics with constant $Q$. See Section \ref{sectionex} 
for additional examples.

For the proof of Theorem \ref{main}, 
we will be concerned with the following equation for a conformal 
metric $\tilde{g} = e^{-2u}g$:
\begin{align}
\label{eqn1}
\sigma_2^{1/2}(g^{-1} A^t_{\tilde{g}}) = f(x) e^{2u}, 
\end{align}
where $f(x) > 0$. We have the following 
formula for the transformation of $A^t$
under a conformal change of metric $\tilde{g} = e^{-2u} g$:
\begin{align}
\label{change1}
A^t_{\tilde{g}} = A^t_g + \nabla^2 u +\frac{1-t}{2}(\Delta u)g
+ du \otimes du - \frac{2-t}{2} |\nabla u|^2 g.
\end{align}
Since $A^t = A^1 + \frac{1-t}{2} tr(A^1) g$,
this formula follows easily from 
the standard formula for the transformation of
the Schouten tensor (see \cite{Jeff1}):
\begin{align}
A^1_{\tilde{g}} = A^1_g + \nabla^2 u + du \otimes du
- \frac{1}{2}| \nabla u|^2 g.
\end{align}
Using (\ref{change1}), we may write (\ref{eqn1}) with 
respect to the background metric $g$
\begin{align}
\label{PDE}
 \sigma_2^{1/2} \left( g^{-1} \Big( \nabla^2 u + \frac{1-t}{2}(\Delta u)g
- \frac{2-t}{2} |\nabla u|^2 g + du \otimes du 
+ A^t_g \Big) \right) = f(x) e^{2u}. 
\end{align}
The choice of the right hand side in (\ref{PDE}) is quite flexible; 
the key requirement is simply that the exponent is a positive multiple 
of $u$.  For negative
exponents we lose the invertibility of the linearized equation 
and some key {\it a priori} estimates; see the proofs of 
Propositions \ref{linvert} and \ref{upper}.

Equation (\ref{PDE}) was considered in our earlier work
 (\cite{GVNegative}) in the context of negative curvature.  
Li and Li (\cite{LiLi2}) used a similar path to prove 
existence of solutions of the 
conformally invariant equation involving more general symmetric 
functions of the eigenvalues, assuming the manifold is locally 
conformally flat.
After completing this paper, we also received the preprint of
Guan, Lin and Wang (\cite{GuanLinWanglcf}),
where they used a similar deformation technique to obtain various 
results in the locally conformally flat setting. 

We will use the continuity method: the assumption of positive 
scalar curvature will allow us to start at some $t = \delta$
very negative. We will then use the conformally 
invariant assumption (\ref{assump}) in Section \ref{C0sec}, 
together with the Harnack inequality 
of \cite{GuanWang1} and \cite{LiLi2} in Section \ref{C1}, 
to prove compactness of the space of solutions. Existence 
of a solution at time $t_0$ and verification of the 
inequalities (\ref{Ricineq}) will be proved in Section~\ref{Existence},
thus completing the proof of Theorem \ref{main}. 
Theorem \ref{Q} will be proved in Section \ref{Pan}, and 
in Section \ref{sectionex} we give many new examples of manifolds 
admitting metrics with constant $Q$-curvature. 
\section{Ellipticity}
\label{ellipticity}
In this section we will discuss the ellipticity 
properties of equation (\ref{PDE}). 
\begin{definition}
Let $(\lambda_1, \dots, \lambda_4) \in \mathbf{R}^4$.
We view the second elementary symmetric function as
a function on $\mathbf{R}^4$:
\begin{align}
\sigma_2(\lambda_1, \dots, \lambda_4) = \sum_{i < j}
\lambda_{i} \lambda_{j},
\end{align}
and we define
\begin{align}
\Gamma_2^+ = \{\sigma_2 > 0\} \cap  \{\sigma_1 > 0\},
\end{align}
where $\sigma_1 = \lambda_1 + \dots + \lambda_4$ denotes the 
trace. 
\end{definition}
For a symmetric linear transformation $A : V \rightarrow V$, where
$V$ is an $n$-dimensional inner product space, the
notation $A \in \Gamma_2^{+}$ will mean that
the eigenvalues of $A$ lie in the corresponding set.
We note that this notation also makes sense for a symmetric tensor on a
Riemannian manifold.
If $A \in \Gamma_2^+$, let
$\sigma_2^{1/2}(A) = \{ \sigma_2(A) \}^{1/2}$.
\begin{definition}
\label{Newtontensor} Let $A : V \rightarrow V$ be 
a symmetric linear transformation,
where $V$ is an $n$-dimensional inner product space.
The {\em{first Newton transformation}} associated with $A$ is
\begin{align}
T_1(A) = \sigma_1(A) \cdot I - A. 
\end{align}
Also, for $t \in  \mathbf{R}$ we define the linear transformation
\begin{align}
\label{Ldef}
L^t(A) = T_{1}(A) + \frac{1-t}{2} \sigma_1(T_1(A)) \cdot I.
\end{align}
\end{definition}
We note that if $A_s : \mathbf{R} \rightarrow
\mbox{Hom}(V,V)$, then
\begin{align}
\label{dsigmak}
\frac{d}{ds} \sigma_2(A_s) = \sum_{i,j} T_{1} (A_s)_{ij} 
\frac{d}{ds}( A_s)_{ij},
\end{align}
that is, the first Newton transformation is what arises from
differentiation of $\sigma_2$. 
\begin{proposition}
\label{ellsumm}
(i)  The set $\Gamma_2^{+}$ is an open convex cone
with vertex at the origin.\\
(ii) If the eigenvalues of $A$ are in $\Gamma_2^{+}$, then 
$T_{1}(A)$ is positive definite. Consequently, for $t \leq 1$, $L^t(A)$ is also 
positive definite.  \\     
(iii)  For symmetric linear transformations $A \in \Gamma_2^+$,
$B \in \Gamma_2^+$, and $s \in [0,1]$,
we have the following inequality
\begin{align}
\label{convexity}
\{ \sigma_2( (1-s)A + sB) \}^{1/2}
\geq (1-s) \{ \sigma_2(A) \}^{1/2}
+ s \{ \sigma_2(B) \}^{1/2}.
\end{align}
\end{proposition}
\begin{proof}
The proof of this proposition is standard, and may be found in
\cite{CNSIII} and \cite{Garding}. 
\end{proof}

For $u \in C^2(M)$, we define 
\begin{align}
\label{changekk}
A^t_u&= A^t_g + \nabla^2 u +\frac{1-t}{2}(\Delta u)g
+ du \otimes du - \frac{2-t}{2} |\nabla u|^2 g.
\end{align}
\begin{proposition}
\label{linvert}
Let $u \in C^2(M)$ be a solution of 
\begin{align}
\label{eqnkk}
\sigma_2^{1/2}(g^{-1} A^t_u ) = f(x) e^{2u}, 
\end{align}
for some $t \leq 1$
with $A_u^t \in \Gamma_2^+$. Then the linearized operator at 
$u$,  $\mathcal{L}^t : C^{2, \alpha}(M) \rightarrow C^{\alpha}(M)$,
is invertible $(0 < \alpha < 1)$. 
\end{proposition}
\begin{proof}
We define 
\begin{align*}
F_t[x, u, \nabla u, \nabla^2 u]
=  \sigma_2 ( g^{-1}A^t_u ) - f(x)^2 e^{4u},
\end{align*}
so that solutions of (\ref{eqnkk}) are zeroes of $F_t$. 
We then suppose that $u \in C^2(M) $ satisfies 
$F_t[x,u, \nabla u, \nabla^2 u] = 0$, with $A^t_u \in \Gamma_2^+$. 
Define $u_s = u + s \varphi$, then 
\begin{align}
\label{linear}
\begin{split}
\mathcal{L}^t(\varphi) &= \frac{d}{ds}F_t[x, u_s, \nabla u_s, \nabla^2 u_s] 
\Big|_{s=0}\\
&= \frac{d}{ds} \Big( \sigma_2 (  g^{-1} A^t_{u_s})\Big) \Big|_{s=0}
- \frac{d}{ds} \left( f^2 e^{4u_s} \right) \Big|_{s=0}.
\end{split}
\end{align}
From (\ref{dsigmak}), we have (using the summation convention)
\begin{align*}
\frac{d}{ds} \Big( \sigma_2 (  g^{-1} A^t_{u} ) \Big) \Big|_{s=0}
= T_{1} ( g^{-1} A^t_{u} )_{ij} \frac{d}{ds} 
\Big( (g^{-1}A^t_{u_s})_{ij} \Big) \Big|_{s=0} 
\end{align*}
We compute
\begin{align*}
  \frac{d}{ds} 
\Big( (g^{-1}A^t_{u_s}) \Big) \Big|_{s=0}   
= g^{-1} \Big( \nabla^2 \varphi  + \frac{1-t}{2}(\Delta \varphi)g
- (2-t)\langle du, d \varphi \rangle g + 
2 du \otimes d \varphi \Big).
\end{align*}
Therefore,
\begin{align}
\notag
\frac{d}{ds} \Big( \sigma_2 (  g^{-1} A^t_{u_s} ) \Big) \Big|_{s=0}
= T_{1} ( g^{-1} A^t_{u} )_{ij}
 \{ g^{-1}  ( 
 & \nabla^2 \varphi   + (1-t)(\Delta \varphi)(g/2)\\
\label{1stterm}
& - (2-t)\langle du, d \varphi \rangle g + 
2 du \otimes d \varphi ) \}_{ij}.
\end{align}
For the second term on the right hand side of (\ref{linear}) we have
\begin{align}
\label{2ndterm}
\frac{d}{ds} \left( f^2 e^{4u_s} \right) \Big|_{s=0}
= 4f^2 e^{4 u} \varphi.
\end{align}
Combining (\ref{1stterm}) and (\ref{2ndterm}), we conclude
\begin{align}
\label{ker}
\mathcal{L}^t(\varphi) = T_{1}( g^{-1} A^t_{u})_{ij}  
\{ g^{-1} ( 
 & \nabla^2 \varphi   + (1-t)(\Delta \varphi)(g/2) ) \}_{ij}
-  4f^2 e^{4 u} \varphi
+ \cdots
\end{align}
where $+ \cdots$ denotes additional terms which are linear in $\nabla \varphi$.
Using the definition of $L^t$ in (\ref{Ldef}), we can rewrite the leading term of
(\ref{ker}) and obtain
\begin{align}
\label{ker2}
\mathcal{L}^t(\varphi) =  L^t(  g^{-1} A^t_{u})_{ij}
(g^{-1} \nabla^2 \varphi)_{ij}  -  4f^2 e^{4 u} \varphi
+ \cdots
\end{align}

For $t \leq 1$, Proposition \ref{ellsumm} implies that 
$L^t(g^{-1} A^t_{u})$ is positive definite, so $\mathcal{L}^t$ is elliptic.  
Since the coefficient of $\varphi$ in the 
zeroth-order term of (\ref{ker2}) is strictly 
negative, the lineariztion is furthermore invertible on the 
stated H\"older spaces (see \cite{GT}). 
\end{proof}
\section{$C^0$ estimate}
\label{C0sec}
Throughout the sequel, $(M,g)$ will be a closed 4-dimensional Riemannian 
manifold with positive scalar curvature.  
Since $R_{g} > 0$, there exists $\delta > - \infty$  so 
that $A^\delta_g$ is positive definite. 
For $t \in [\delta, 1]$, consider the path of equations
\begin{align}
\label{path}
\sigma_2^{1/2}( g^{-1}A^t_{u_t} ) = f(x) e^{2u_t},
\end{align}
where $f(x) =\sigma_2^{1/2}(g^{-1} A^\delta_g ) > 0$.
Note that $u \equiv 0$ is a solution of (\ref{path}) for $t=\delta$. 
\begin{proposition} 
\label{upper}
Let $u_t \in C^2(M)$ be a solution of (\ref{path})
for some $\delta \leq t \leq 1$.
Then $u_t \leq \overline{\delta}$, where $\overline{\delta}$ depends 
only upon $g$.
\end{proposition} 
\begin{proof}
From Newton's inequality $\frac{4}{\sqrt{6}} \sigma_2^{1/2} 
\leq  \sigma_1$, so 
\begin{align}
\frac{4}{\sqrt{6}} f(x) e^{2u_t} \leq \sigma_1( g^{-1}A^t_{u_t} ).
\end{align}
Let $p$ be a maximum of $u_t$, then the gradient terms 
vanish at $p$, and $\Delta u_t \leq 0$, so by (\ref{change1})
\begin{align*}
\begin{split}
\frac{4}{\sqrt{6}} f(p) e^{2u_t(p)} & \leq \sigma_1( g^{-1}A^t_{u_t} )(p) \\
&= \sigma_1(g^{-1} A_g^t) + (3-2t) \Delta u_t \\
&\leq \sigma_1(g^{-1}A_g^t).
\end{split}
\end{align*}
Since $t \geq \delta$, this implies $u_t \leq \overline{\delta}$. 
\end{proof}
\begin{proposition}
\label{C0}
Assume that for some $\delta \leq t \leq 1$, 
\begin{align} 
\mathcal{F}_2([g]) + \frac{1}{6} (1-t) (2-t) (Y[g])^2 = \lambda_t > 0.
\end{align}
If $u_t \in C^2(M)$ is a solution of (\ref{path})
satisfying $\Vert \nabla u_t \Vert_{L^{\infty}} < C_1$,
then $ u_t > \underline{\delta}$, 
where $\underline{\delta}$ depends only upon $g, C_1,$ and $\log \lambda_t $. 
\end{proposition}
\begin{proof}
Using Lemma 24 in \cite{Jeff1}, we have 
\begin{align*}
\sigma_2(A^t) &= \sigma_2 \Big( A^1 + \frac{1-t}{2} \sigma_1(A^1) g \Big)\\
& =  \sigma_2 ( A^1) + 3 \frac{1-t}{2} \sigma_1(A^1)^2 
+ 6  \Big( \frac{1-t}{2} \sigma_1(A^1) \Big)^2\\
& =  \sigma_2 ( A^1) +\frac{3}{2} (1-t)(2-t) \sigma_1(A^1)^2.
\end{align*}
Letting $\tilde{g} = e^{-2 u_t} g$, 
\begin{align*}
e^{4u_t} f^2 =\sigma_2(g^{-1}A^t_{u_t}) & = \sigma_2(g^{-1}A^1_{u_t}) 
+ \frac{3}{2}(1-t)(2-t) \big (\sigma_1(g^{-1}A^1_{u_t}) \big)^2\\
&=  e^{-4 {u}_t} \Big( \sigma_2( \tilde{g}^{-1}A^1_{u_t}) +  
 \frac{1}{24}(1-t)(2-t) R_{\tilde{g}}^2 \Big).
\end{align*}
Integrating this, we obtain
\begin{align*}
C' \int_M e^{4u_t} dvol_{g} &\geq 
\int_M f^2 e^{4u_t} dvol_{g} \\
& = \int_M 
\sigma_2( \tilde{g}^{-1}A^1_{u_t}) e^{-4 u_t} dvol_{g}
+  \frac{1}{24}(1-t)(2-t)  \int_M   R_{\tilde{g}}^2  e^{-4 u_t} dvol_{g}\\
& = \int_M \sigma_2( \tilde{g}^{-1}A^1_{\tilde{g}}) dvol_{\tilde{g}}
+  \frac{1}{24}(1-t)(2-t)  \int_M   R_{\tilde{g}}^2  dvol_{\tilde{g}},
\end{align*}
where $C' > 0$ is chosen so that $f^2 \leq C'$.
\begin{lemma} For any metric $g' \in [g]$, we have 
\begin{align}
\int_M R_{g'}^2 dvol_{g'} \geq (Y[g])^2.
\end{align}
\end{lemma}
\begin{proof}
From H\"older's inequality, 
\begin{align}
\label{yadda}
\int_M R_{g'} dvol_{g'} \leq \left\{ \int_M R_{g'}^2 dvol_{g'} \right\}^{1/2}
\cdot \left\{ Vol(g') \right\}^{1/2}. 
\end{align}
Since $g$ has positive scalar curvature, $Y[g] > 0$, so the 
left hand side of (\ref{yadda}) must be positive. 
We then obtain
\begin{align*}
(Y[g])^2 \leq  \left( Vol(g')^{-1/2} \int_M R_{g'} 
dvol_{g'} \right)^2 \leq \int_M R_{g'}^2 dvol_{g'}. 
\end{align*}
\end{proof}
Using the lemma, and the conformal invariance of 
$\mathcal{F}_2$, we obtain
\begin{align}
\label{google}
C' \int_M e^{4u_t} dvol_{g}
&\geq  \frac{1}{4} \mathcal{F}_2([g]) 
+  \frac{1}{24}(1-t)(2-t) (Y[g])^2 = \frac{1}{4}\lambda_t > 0.
\end{align}
This implies
\begin{align}
\label{supbelow}
\max u_t \geq \frac{1}{4} \log \lambda_t - C(g).
\end{align}
The assumption $|\nabla u_t| < C_1$ implies the Harnack inequality
\begin{align}
\label{Harnack}
\max u_t \leq \min u_t + C(C_1,g),
\end{align}
by simply integrating along a geodesic connecting points at which $u_t$
attains its maximum and minimum.
Combining (\ref{supbelow}) and (\ref{Harnack}) we obtain
$$
\min u_t \geq \frac{1}{4} \log\lambda_t - C.
$$
\end{proof}
\section{Harnack inequality}
\label{C1}
We next have the following $C^1$ estimate for solutions of the equation 
(\ref{PDE}).
\begin{proposition}
\label{C1estimate}
 Let $u_t$ be a $C^3$ solution of (\ref{path}) for some $\delta \leq t \leq 1$, 
satisfying $u_t < \overline{\delta}$.
Then $\Vert \nabla u_t \Vert_{L^{\infty}} < C_1$, where $C_1$ 
depends only upon $\overline{\delta}$ and $g$. 
\end{proposition}
\begin{remark}
A Harnack inequality was proved for the conformally invariant 
equation for $t=1$ in \cite{GuanWang1}, and then extended to 
$t < 1$ in \cite{LiLi2}. More specifically, in \cite{LiLi2} was
considered
the equation
\begin{align}
\sigma_k^{1/k}( s A^1 + (1-s) \sigma_1(A^1) g) = f(x) e^{-2u}.
\end{align}
The left hand side is just a reparametrization of $A^t$,
but (\ref{path}) has a different 
right hand side, so the Harnack inequality now depends on the 
$sup$. The differences are minor, but for convenience, 
we present an outline of the proof here, and also provide a
simple direct proof which works for $t < 1$. 
\end{remark}
\begin{proof}
Consider the function $h = |\nabla u|^2$ (we will omit 
the subscript on $u_t$). 
Since $M$ is compact, and $h$ is 
continuous, we suppose the maximum of $h$ occurs and a point 
$p \in N$. Take a normal coordinate system 
$(x^1, \dots, x^n)$ at $p$, then 
$g_{ij}(p) = \delta_{ij}$, and $\Gamma^i_{jk}(p) = 0$, 
where $g = g_{ij} dx^i dx^j$, and $\Gamma^i_{jk}$
is the Christoffel symbol (see \cite{Besse}). 

 Locally, we may write $h$ as 
\begin{align}
h =g^{lm}u_{l}u_{m}.
\end{align}
In a neighborhood of $p$, differentiating $h$ in the $x^i$
direction we have
\begin{align}
\partial_i h = h_i = \partial_i ( g^{lm}u_l u_m)
\label{goog} = \partial_i ( g^{lm})u_l u_m  
 + 2 g^{lm} \partial_i(u_l) u_m   
\end{align}
Since in a normal coordinate system, the first 
derivatives of the metric vanish at $p$, 
and since $p$ is a maximum for $h$, 
evaluating (\ref{goog}) at $p$, we have
\begin{align}
\label{hip}
u_{li}u_l = 0. 
\end{align}
Next we differentiate (\ref{goog}) in the $x^j$ direction. 
Since $p$ is a maximum, $\partial_j \partial_i h = h_{ij}$ is negative 
semidefinite, and we get (at $p$)
\begin{align}
\label{zonk}
0 \gg h_{ij} &=  \frac{1}{2}\partial_j \partial_i g^{lm} u_l u_m
 + u_{lij}u_l  + u_{li}u_{lj}. 
\end{align} 
We recall from Section \ref{ellipticity} that
\begin{align}
\label{Lagain}
L^t_{ij} = {T}_{ij} + \frac{1-t}{2} \sum_l {T}_{ll} \delta_{ij}, 
\end{align}
is positive definite, where ${T}_{ij}$
means $(T_{1}(g^{-1}A^t_{u}))_{ij}$.
We sum with (\ref{zonk}) with $L^t_{ij}$ to obtain
the inequality
\begin{align}
\label{inequality1}
0 \geq \frac{1}{2} L^t_{ij} \partial_i \partial_j
g^{lm}u_l u_m + L^t_{ij} u_{lij}u_l +  L^t_{ij} u_{li} u_{lj}. 
\end{align}

We next differentiate equation (\ref{path}) in order to replace 
the $u_{lij}$ term with lower order terms. With respect to 
our local coordinate system, from (\ref{changekk}) we have 
\begin{align}
\notag
(A^t_u)_{ij} = &(A^t_g)_{ij} + u_{ij} - u_r \Gamma^r_{ij} 
+\frac{1-t}{2}\sum_k( u_{kk} - u_r \Gamma^r_{kk})g_{ij}
+ u_i u_j \\
\label{ubarlocal}
&- \frac{2-t}{2} (g^{r_1 r_2} u_{r_1} u_{r_2} )g_{ij}.
\end{align} 
At the point $p$, this simplifies to 
\begin{align}
\label{uijeqn}
(A^t_u)_{ij} =  (A^t_g)_{ij}
+ u_{ij}  +\frac{1-t}{2}\sum_k(u_{kk})g_{ij}
+ u_i u_j - \frac{2-t}{2} (|\nabla u|^2 )\delta_{ij}.
\end{align}

Next we take $m$ with $1 \leq m \leq n$, and differentiate 
(\ref{path}) with respect to $x^m$ in our local coordinate system:
\begin{align}
\label{firstpartial}
\partial_m \left\{ \sigma_2 \big( g^{lj}( A^t_u)_{ij}  \big) \right\} = 
\partial_m (f(x)^2 e^{4u}).
\end{align}
Differentiating and evaluating at $p$, we obtain
\begin{align}
\begin{split}
\label{fullfirstder}
&{T}_{ij} \Big( 
\partial_m (A^t_g)_{ij}
+ u_{ijm} - u_r \partial_m  \Gamma^r_{ij} 
+\frac{1-t}{2} \sum_k (u_{kkm} - u_{r} \partial_m \Gamma^r_{kk})\delta_{ij} 
+ 2u_{im} u_j \Big)\\
& = (\partial_m f^2) e^{4u}
+ 4 f^2 e^{4u}u_m. 
\end{split}
\end{align}
Note that the third order terms in the above expression are 
\begin{align*}
&{T}_{ij} \Big( u_{ijm} +\frac{1-t}{2}\sum_k u_{kkm}\delta_{ij} \Big)
= L^t_{ij} u_{ijm}.    
\end{align*}
Next we sum (\ref{fullfirstder}) with $u_m$, using (\ref{hip})
we have the following formula
\begin{align}
\begin{split}
\label{firstderterm}
& L^t_{ij} u_m u_{ijm} + 
{T}_{ij} \Big( 
 u_m \partial_m (A^t_g)_{ij} - u_m u_r \partial_m  \Gamma^r_{ij} 
-\frac{1-t}{2}\sum_k( u_{r} u_m \partial_m \Gamma^r_{kk})\delta_{ij} \Big)\\
& \hspace{30mm} = u_m (\partial_m f^2) e^{4u} 
+ 4 f^2 e^{4u} | \nabla u|^2. 
\end{split}
\end{align}
Substituting (\ref{firstderterm}) into 
(\ref{inequality1}), we arrive at the inequality
\begin{align*}
0 &\geq \frac{1}{2}  L^t_{ij} \partial_i \partial_j
 g^{lm}u_l u_m 
+ {T}_{ij} \Big( -
u_m \partial_m (A^t_g)_{ij} + u_m u_r \partial_m  \Gamma^r_{ij} 
+\frac{1-t}{2}\sum_k( u_{r} u_m \partial_m \Gamma^r_{kk})\delta_{ij}
 \Big)\\
&+u_m (\partial_m f^2) e^{4u} + L^t_{ij} u_{li} u_{lj}.
\end{align*}
Using (\ref{Lagain}) and Lemma 2 in \cite{Jeff4}, we obtain
\begin{align}
\notag
0 &\geq {T}_{ij} \Big( \frac{1-t}{2}  \sum_k R_{klkm} {u_l u_m} \delta_{ij}
+ R_{iljm} {u_l u_m} - {u_m} \partial_m (A^t_g)_{ij} \Big)\\
\label{mess2}
& \ \ \ \ \ \ + u_m (\partial_m f^2) e^{4u}
+ {T}_{ij}u_{li}u_{lj} + \frac{1-t}{2} \sum_l {T}_{ll}u_{ij}u_{ij},
\end{align}
where $R_{iljm}$ are the components of the Riemann
curvature tensor of $g$.

\begin{lemma}
\label{extraterm} 
There exists a constant $\beta > 0$ such that for $t \in [\delta,1]$,  
\begin{align}
{T}_{ij}u_{li}u_{lj} + \frac{1-t}{2} \sum_l{T}_{ll}u_{ij}u_{ij}
\geq \beta \sum_l {T}_{ll} |\nabla u|^4
\end{align}
\end{lemma}

\begin{remark}
This was proved in \cite{LiLi2}, using the result in \cite{GuanWang1}.
\end{remark}

Using the lemma, we have 
\begin{align}
\notag
0 \geq {T}_{ij} & \Big( \frac{1-t}{2}  \sum_k R_{klkm} {u_l u_m} \delta_{ij}
+ R_{iljm} {u_l u_m} - {u_m} \partial_m (A^t_g)_{ij} \Big)\\
&+ u_m (\partial_m f^2) e^{4u} 
+ \beta \sum_l {T}_{ll}| \nabla u |^4.
\end{align}
Since we are assuming $u$ is bounded above, the $|\nabla u|^4$ 
term dominates, and the proof proceeds as in \cite{GuanWang1} or 
\cite{LiLi2}. 
\end{proof}

\begin{remark}

For convenience, we would like to present here a simplified proof of Lemma \ref{extraterm} which works
for $t < 1$.  Although the argument breaks down as $t \rightarrow 1$,
it covers the case $t_0 = 0$, and therefore suffices for    
proving Theorem \ref{Q}.

To begin, we claim that if $\beta_t' > 0 $ is sufficiently small, 
then for at least one $i_0$, $u_{i_0 i_0} \geq \beta_t' |\nabla u|^2$.  
If not,  then $u_{ii} < \beta_t' |\nabla u|^2$
for $i = 1 \dots 4$.
Since $\Gamma_2 \subset \{ \sigma_1 > 0 \}$, 
\begin{align*}
 \Delta u + 2 (1-t)(\Delta u)
- 2(2-t) |\nabla u|^2 + |\nabla u|^2 
+ \sigma_1(A^t_g)  > 0.
\end{align*}
Without loss of generality, we may assume that 
$\sigma_1(A^t_g) \leq \epsilon |\nabla u|^2$, we then have 
\begin{align*}
\Big( 1 +  2 (1-t) \Big) \Delta u 
+ \Big( 1 - 2(2-t) \Big) |\nabla u |^2 + \epsilon |\nabla
u|^2 > 0.
\end{align*}
From the assumption, $ \Delta u = 
\sum u_{ii} < 4 \beta_t' |\nabla u|^2$, so 
\begin{align*}
4 \beta_t' \Big( 1 +  2(1-t) \Big) |\nabla u|^2 
+ \Big( 1 - 2 (2-t) \Big) |\nabla u |^2 + 
\epsilon |\nabla u|^2 > 0, 
\end{align*}
which is a contradiction for $\epsilon$ and $\beta_t'$ sufficiently small.
We then have 
\begin{align}
{T}_{ij}u_{li}u_{lj} + \frac{1-t}{2} \sum_l {T}_{ll}u_{ij}u_{ij}
\geq  \frac{1-t}{2} \sum_l {T}_{ll}u_{i_0 i_0}^2
\geq  \frac{1-t}{2} (\beta_t')^2 \sum_l {T}_{ll} |\nabla u|^4,
\end{align}
so choose $\beta_t =  \frac{1-t}{2} (\beta_t')^2$.  This completes the proof.
\end{remark}

\section{Proof of Theorem \ref{main}}
\label{Existence}
\begin{proposition}
\label{C2estimate}
 Let $u_t$ be a $C^{4}$ solution of (\ref{path}) for some $\delta \leq 
t \leq 1$ satisfying $\underline{\delta} < u_t < \overline{\delta}$,
and $\Vert \nabla u_t \Vert_{L^{\infty}} < C_1$.
Then for $0 < \alpha < 1$,  $\Vert u_t \Vert_{C^{2, \alpha}} \leq C_2$,
where $C_2$ 
depends only upon $\underline{\delta},\overline{\delta}, 
 C_1$, and $ g$. 
\end{proposition}
\begin{proof}
The $C^2$ estimate follows from the global 
estimates in \cite{GVNegative}, or the local estimates  
\cite{GuanWang1} and \cite{LiLi2}. 
We remark that the main fact used in deriving these estimates is 
that $\sigma_2^{1/2}(A^t)$ is a {\em{concave}} function of the second 
derivative variables, which follows easily from the
inequality (\ref{convexity}). 
Since $f(x) > 0$, the $C^2$ estimate implies uniform 
ellipticity, 
and the $C^{2, \alpha}$ estimate then follows 
from the work of \cite{Krylov} and \cite{Evans} on 
concave, uniformly elliptic equations. 
\end{proof}

 To finish the proof of Theorem \ref{main}, we use the
continuity method. Recall that we are considering the 
1-parameter family of equations, for $t \in [\delta, t_0]$, 
\begin{align}
\label{path2}
\sigma_2^{1/2}( g^{-1}A^t_{u_t} ) = f(x) e^{2u_t},
\end{align}
with $f(x) =\sigma_2^{1/2}(g^{-1} A^\delta_g ) > 0$,
and $\delta$ was chosen so that $A^\delta_g$ is positive definite.
We define
\begin{align*}
\mathcal{S} = \{ t \in [ \delta, t_0]
 \ | \ \exists \mbox{ a solution } u_t  \in C^{2, \alpha}(M) 
\mbox{ of } (\ref{path2})
\mbox{ with } A^t_{u_t} \in \Gamma_2^+ \}.
\end{align*}
The function $f(x)$ was chosen so that  
$u \equiv 0$ is a solution at $t = \delta$. 
Since $A^{\delta}_g$ is positive definite, 
and the positive cone is clearly contained in  $\Gamma_2^+$,
$\mathcal{S}$ is nonempty. Let $t \in \mathcal{S}$, and $u_t$
be any solution. From Proposition \ref{linvert}, the linearized 
operator at $u_t$, $\mathcal{L}^t: C^{2,\alpha}(M) \rightarrow 
C^{\alpha}(M)$, is invertible. The implicit function 
theorem (see \cite{GT}) implies that 
$\mathcal{S}$ is open. Note that since $f \in C^{\infty}(M)$, 
it follows from classical elliptic regularity theory 
that $u_t \in C^{\infty}(M)$. 
Proposition \ref{upper} implies a
uniform upper bound on solutions $u_t$ (independent 
of $t$). We may then apply Proposition \ref{C1estimate} 
to obtain a uniform gradient bound, and Lemma \ref{C0} 
then implies a uniform lower bound on $u_t$. 
Proposition \ref{C2estimate} then implies that 
$\mathcal{S}$ is closed, therefore $ \mathcal{S} = [\delta,t_0]$. 
The metric $\tilde{g} = e^{-2u_{t_0}} g$ then satisfies
$\sigma_2( A^{t_0}_{\tilde{g}} ) > 0$ and $ R_{\tilde{g}} > 0$. 

  We next verify the inequalities (\ref{Ricineq}). 
We decompose $A^t$ into its trace-free and pure-trace components,
\newcommand{\oca}{\overset{\circ}{A^t}}
\begin{align}
\label{decomp}
\begin{split}
A^t &= A^t - \frac{1}{n}\sigma_1(A^t)g + \frac{1}{n}\sigma_1(A^t)g \\
& \equiv \oca + \frac{1}{n}\sigma_1(A^t)g.
\end{split}
\end{align}
We now associate to $A^t$ the
symmetric transformation $\widehat{A^t}$, defined by
\begin{align}
\label{Ahat}
\widehat{A^t} \equiv -\oca + \frac{1}{n}\sigma_1(A^t)g.
\end{align}
That is, $\widehat{A^t}$ is the (unique) symmetric transformation 
which has the same pure-trace component as $A$, 
but the opposite trace-free component.  
\begin{lemma}
\label{Ahatprop}The tensors  $\widehat{A^t}$ and $A^t$ satisfy
the equalities
\begin{align}
\label{i}
\sigma_1(\widehat{A^t}) = \sigma_1(A^t), \\
\label{ii}
\sigma_2(\widehat{A^t}) = \sigma_2(A^t).
\end{align}
\end{lemma}
\begin{proof}
The proof of $(\ref{i})$ is immediate from the definition of $\widehat{A^t}$.  
To prove $(\ref{ii})$,
we use the identity
\begin{align*}
\sigma_2(A^t) &= -\frac{1}{2}|A^t|^2 + \frac{1}{2}\sigma_1(A^t)^2 \\ 
&= -\frac{1}{2}|\oca + \frac{1}{n}\sigma_1(A^t)g|^2 
+ \frac{1}{2}\sigma_1(A^t)^2.
\end{align*}
Since the decomposition (\ref{decomp}) is orthogonal with respect to the
norm $|\cdot|^2$, we conclude
\begin{align*}
\sigma_2(A^t)
= -\frac{1}{2}|-\widehat{A^t}|^2 + \frac{1}{2}\sigma_1(A^t)^2 
= \sigma_2(\widehat{A^t}).
\end{align*}
\end{proof}

Combining Lemma \ref{Ahatprop} and Proposition \ref{ellsumm} we have

\begin{proposition}
\label{otherineq}
If the eigenvalues of $A^t$ are in $\Gamma_2^{+}$, then 
\begin{align}
\label{i1}
- &A^t + \sigma_1(A^t)g > 0, \mbox{ and} \\
\label{i2}
&A^t + \frac{n-2}{n}\sigma_1(A^t)g > 0.
\end{align}
\end{proposition}
\begin{proof}
The tensor in (\ref{i1}) is simply the first Newton 
transformation of $A^t$, which is positive definite 
by Proposition \ref{ellsumm}.
By Lemma \ref{Ahatprop}, $\sigma_2(\widehat{A^t}) > 0$
and $\sigma_1(\widehat{A^t}) > 0$.  Thus, the eigenvalues of $\widehat{A^t}$
are also in $\Gamma_2^{+}$.
By Proposition \ref{ellsumm}, 
the first Newton tranform $T_1(\widehat{A^t})$ is positive definite.  
By definition,
\begin{align*}
T_1(\widehat{A^t}) &= -\widehat{A^t} + \sigma_1(\widehat{A^t})g \\
&= -\left( -\oca + \frac{1}{n}\sigma_1(A^t)g\right) + \sigma_1(A^t)g \\
&= A^t + \frac{(n-2)}{n}\sigma_1(A^t)g. 
\end{align*}
\end{proof}
When expressed in terms of $Ric$ and $n=4$, (\ref{i1}) and (\ref{i2}) 
are exactly (\ref{Ricineq}).
\section{Proof of Theorem \ref{Q}}
\label{Pan}
The assumption (\ref{assump3}) corresponds to $t_0=0$ in 
Theorem \ref{main}. From (\ref{Ricineq}) we find a conformal 
metric (which for simplicity, we again denote by $g$) 
with $ Ric_g < R_g g$. Theorem \ref{Q} then follows 
from  
\begin{proposition}
\label{ricr}
If  $ Ric_g \leq R_g g$, then $P \geq 0$, and $Ker P = \{ constants \}$. 
\end{proposition}
\begin{proof}
We again recall that the Paneitz operator is defined by 
\begin{align}
P \phi   = \Delta^2 \phi + \delta \Big( \frac{2}{3} R_g g - 2 Ric_g \Big)
d\phi. 
\end{align}
Integrating by parts, 
\begin{align}
\label{intp}
\langle P \phi, \phi \rangle_{L^2}
= \int_M \left( (\Delta \phi)^2 + \frac{2}{3} R_g |\nabla \phi|^2
- 2 Ric_g( \nabla \phi, \nabla \phi) \right) dvol_g. 
\end{align}
From the Bochner formula, 
\begin{align}
\label{Bochner}
0 = \int \left( | \nabla^2 \phi|^2 + Ric_g( \nabla \phi, \nabla \phi) 
- ( \Delta \phi)^2 \right) dvol_g.
\end{align}
Substituting (\ref{Bochner}) into (\ref{intp}), we have
\begin{align*}
\langle P \phi, \phi \rangle_{L^2}
&= \int_M \left(-\frac{1}{3} (\Delta \phi)^2 +\frac{4}{3} (\Delta \phi)^2
+ \frac{2}{3} R_g |\nabla \phi|^2 - 2 Ric_g( \nabla \phi, \nabla \phi)
\right) dvol_g\\
& =  \int_M \left( -\frac{1}{3} (\Delta \phi)^2 +\frac{4}{3} |\nabla^2 \phi|^2
+ \frac{2}{3} R_g |\nabla \phi|^2 - \frac{2}{3} 
Ric_g( \nabla \phi, \nabla \phi) \right) dvol_g\\
& =  \int_M \left( \frac{4}{3} | \overset{\circ}{\nabla^2} \phi|^2
+ \frac{2}{3} ( R_g g - Ric_g)( \nabla \phi, \nabla \phi) \right) dvol_g\\
& \geq \int_M \frac{4}{3} | \overset{\circ}{\nabla^2} \phi|^2 dvol_g,
\end{align*}
where $ \overset{\circ}{\nabla^2} \phi = \nabla^2 \phi
- (1/4) (\Delta \phi) g$. Consequently, $P \geq 0$. 
Assume by contradiction that $P \phi = 0$, and $\phi$ is not 
constant. From the above, we conclude
that $ \overset{\circ}{\nabla^2} \phi \equiv 0$. 
By \cite[Theorem A]{Obata2}, $g$ is homothetic
to $S^4$. We then have 
\begin{align*}
0 = \langle P \phi, \phi \rangle_{L^2}
& =  \int_M \frac{4}{3} | \overset{\circ}{\nabla^2} \phi|^2 dvol_g
+ \frac{1}{2} R_g \int_M | \nabla \phi|^2 dvol_g,
\end{align*}
and therefore $\phi = constant$. 
\end{proof}
\begin{remark}
In \cite{ES}, the nonnegativity of the Paneitz operator 
was shown assuming $R g - \lambda Ric \geq 0$ for 
$\lambda \in (1,3]$. The above proposition extends this 
to the endpoint $\lambda = 1$. 
\end{remark}
\begin{corollary} 
\label{ricg}
If  $ Ric_g \geq 0$, then $P \geq 0$, and $Ker P = \{ constants \}$. 
\end{corollary}
\begin{proof}
Clearly, $Ric_g \geq 0$ implies that $Ric_g \leq R_g g$, so this 
follows directly from Proposition \ref{ricr}. 
\end{proof}
\begin{remark}
The construction in \cite{ShaYang2} yields metrics with 
positive Ricci curvature on the connect sums 
$k(S\sp 2\times S\sp 2)$, 
$k \mathbb{CP}^2 \#\overline{ \mathbb{CP}}{}\sp 2)$, 
and $(k+l) \mathbb{CP}^2 \#k\overline{ \mathbb{CP}}{} \sp 2$.
Consequently, from Corollary \ref{ricr}, and the results in 
\cite{ChangYang}, these manifolds admit metrics with $Q = constant$. 
\end{remark}
\section{Examples}
\label{sectionex}

The following theorem will allow us to give many examples 
of metrics satisfying the conditions of Theorem \ref{Q}.
\begin{theorem}
\label{examples}
Let $(M,g)$ satisfy $\int Q_g dvol_g \geq 0$.
If $Y[g] > 4 \sqrt{3k} \pi, \ k < 8$,
then the manifold $N = M \# k(S^1 \times S^3)$ admits a metric $\tilde{g}$ 
satisfying (\ref{assump3}).
If $Y[g] > 8 \sqrt{3} \pi$, 
then the manifold $N = M \# l( \mathbb{RP}^4 )$ admits a metric $\tilde{g}$ 
satisfying (\ref{assump3}) for $l < 9$. 
Consequently, these manifolds $N$ admit metrics with $Q = constant.$
\end{theorem}
\begin{proof}
From the assumption on $\int Q dvol $ and the 
Chern-Gauss-Bonnet formula, we have
\begin{align}
\label{Weyless} 
\int_M |W_g|^2 dvol_g \leq 8 \pi^2 \chi(M)
\end{align}
From \cite[Proposition 4.1]{Kobayashi2}, given a point 
$p \in M$, and $\epsilon > 0$, $M$ admits a metric $g'$ so that 
$g'$ is locally conformally flat in a neighborhood 
of $p$, and
\begin{align}
\begin{split}
\label{stuff}
&|Y[g'] - Y[g]| < \epsilon,\\
& \Big| \int_M |W_g|^2 dvol_g - \int_M |W_{g'}|^2 dvol_{g'} \Big| 
< \epsilon.
\end{split}
\end{align} 
We next put a metric on the connect sum using the
technique in \cite{Kobayashi}. Since $g'$ is
locally conformally flat near $p$, there is a conformal
factor on $M - \{p\}$ which makes the metric look cylindrical 
near $p$. 

For the first case, since $S^1 \times S^3$ is locally 
conformally flat, for any  $p' \in S^1 \times S^3$ there is a 
conformal factor on $S^1 \times S^3 - \{p'\}$ which makes the 
metric look cylindrical 
near $p'$. Therefore one can put a metric on $N$ by 
identifying the cylindrical regions together along their
boundaries.  From the construction in \cite{Kobayashi}, there are sequences 
of locally conformally flat metrics on $k(S^1 \times S^3)$ whose 
Yamabe invariants approach 
$\sigma(k(S^1 \times S^3)) = \sigma(S^4) = 8 \sqrt{6} \pi 
> Y[g']$, where 
$\sigma$ denotes the diffeomorphism Yamabe invariant, so we choose 
a locally conformally flat metric $g_1$ on $k(S^1 \times S^3)$ 
satisfying $Y[g_1] \geq 8 \sqrt{6} \pi - \epsilon$. 
We have $min\{ Y[g_1], Y[g'] \} = Y[g']$, 
so following the proof \cite[Theorem 2]{Kobayashi},
by changing the length of the cylindrical region, 
one can put a metric $\tilde{g}$ on the connect sum 
$N = M \# k(S^1 \times S^3)$ with $ Y[ \tilde{g}]  > Y[g'] - \epsilon$.
Clearly we also have 
\begin{align*}
&\Big| \int_M |W_g|^2 dvol_g - \int_N |W_{\tilde{g}}|^2
dvol_{\tilde{g}} \Big| 
< \epsilon,
\end{align*}
which along with (\ref{Weyless}) implies
\begin{align}
\label{Weyltilde}
\int_N |W_{\tilde{g}}|^2 dvol_{\tilde{g}} \leq 8\pi^2\chi(M) + \epsilon.
\end{align}

We next verify that, for appropriate $\epsilon$, 
the metric $\tilde{g}$ satisfies the condition (\ref{assump3}).
To see this, 
write $(Y[g])^2 = 48 k \pi^2 + 3 \delta$, with $\delta > 0$,
and noting that $\chi(N) = \chi(M) - 2k$ we have 
\begin{align*}
2 \int_N Q_{\tilde{g}} dvol_{\tilde{g}} + \frac{1}{3} ( Y[\tilde{g}])^2 &= 
8 \pi^2 \chi(N) - \int_N |W_{\tilde{g}}|^2 dvol_{\tilde{g}}
+ \frac{1}{3} ( Y[ \tilde{g}])^2\\
&\geq  8 \pi^2 \chi(N) -  8 \pi^2 \chi(M) 
+  \frac{1}{3} ( Y[g])^2 - C \epsilon\\
& = -16k \pi^2 + \frac{1}{3} ( Y[g])^2 - C \epsilon \\ 
&=  \delta - C \epsilon > 0, 
\end{align*}
for $\epsilon$ sufficiently small. 

 For the second case, since $\mathbb{RP}^4$ is locally 
conformally flat, we do exactly the same gluing as before. 
Again we use \cite[Proposition 4.1]{Kobayashi2} to 
find a metric $g'$ on $M$ satisfying (\ref{stuff}). 
We fix the standard metric $g_0$ on $\mathbb{RP}^4$,
which is locally conformally flat. 
Since $Y([g']) > 8 \sqrt{3} \pi = Y[g_0]$
we have $min \{ Y[g_0], Y[g'] \} = Y[g_0]$, so from the construction 
in \cite{Kobayashi}, we can put a metric  $\tilde{g}$ on the connect sum 
$N = M \# l(\mathbb{RP}^4)$ with $ Y[ \tilde{g}]  > 8 \sqrt{3} \pi 
- \epsilon$, and which also satisfies (\ref{Weyltilde}). 

Write $(Y[g])^2 = 3( 64 \pi^2 + \delta)$, with $\delta > 0$,
and noting that $\chi(N) = \chi(M) -l$,
\begin{align*}
2 \int_N Q_{\tilde{g}} dvol_{\tilde{g}} + \frac{1}{3} ( Y[\tilde{g}])^2 &=
8 \pi^2 \chi(N) - \int_N |W_{\tilde{g}}|^2 dvol_{\tilde{g}}
+ \frac{1}{3} ( Y[ \tilde{g}])^2\\
&\geq  8 \pi^2 \chi(N) -  8 \pi^2 \chi(M) 
+  \frac{1}{3} ( Y[g])^2 - C \epsilon\\
& = -8l \pi^2 + 64 \pi^2 + \delta - C \epsilon > 0, 
\end{align*}
for $\epsilon$ sufficiently small and $l < 9.$
\end{proof}
We next write down some specific examples of $(M,g)$ satisfying the
assumptions of Theorem \ref{examples}.
We will use the fact that if $(M,g)$ is a positive 
Einstein manifold, then
 $\int Q_g dvol_g > 0$ and the Yamabe invariant
is attained by $g$. 

(1) $M = S^2 \times S^2$ with the product metric,
$Y[g] = 16 \pi >  4 \sqrt{3k} \pi$ for $k < 6$, so we have 
\begin{align}
N =  (S^2 \times S^2) \# k(S^1 \times S^3), \ k \leq 5.
\end{align}

(2) $M =  \mathbb{CP}^2$ with the Fubini-Study metric, 
$Y[g] = 12 \sqrt{2} \pi$ (see \cite{Lebrun}).
Since $12 \sqrt{2} \pi > 4 \sqrt{3k} \pi$ for 
$k < 6$, this yields the examples
\begin{align}
N =  \mathbb{CP}^2 \# k(S^1 \times S^3), \ k \leq 5.
\end{align}

(3) Again, we take  $M =  \mathbb{CP}^2$. We have  
$12 \sqrt{2} \pi > 8 \sqrt{3} \pi$, so from the
second statement in Theorem \ref{examples} we have 
\begin{align}
N =  \mathbb{CP}^2 \# k( \mathbb{RP}^4 ), \ k \leq 8.
\end{align}

(4) $M = \mathbb{CP}^2 \# l \overline{ \mathbb{CP}}{}^2$,
$3 \leq l \leq 8$, $M$ admits K\"ahler-Einstein metrics
satisfying $Y[g] = 4 \pi \sqrt{2(9-l)}$ (see \cite{Gurskyweyl}, 
\cite{Lebrun}). 
Since $4 \pi \sqrt{2(9-l)} > 4 \sqrt{3} \pi$ for $l < 8$, we have
the examples
\begin{align}
N = \mathbb{CP}^2\# l \overline{ \mathbb{CP}}{}^2
\#( S^1 \times S^3), \ 3 \leq l \leq 7.
\end{align}

(5) $N = k(S^1 \times S^3) \#  l(\mathbb{RP}^4), \ 2k +l \leq 9$.
We do not need Theorem \ref{examples} for this example, we
argue directly. By the construction in \cite{Kobayashi}, these 
manifolds admit locally conformally flat metrics $\tilde{g}$ with
$Y[ \tilde{g}] \approx Y[ \mathbb{RP}^4, g_0] = 8 \sqrt{3} \pi$,
We have $\chi(N) = -2k -l +2$, so the assumption (\ref{assump3}) 
is that 
\begin{align*}
0 <  8 \pi^2 (-2k - l + 2) + \frac{1}{3}(Y[g])^2 
\approx  8 \pi^2 (-2k - l + 2) + 64 \pi^2,
\end{align*}
which is satisfied for $2k + l < 10$.

The above examples are all summing with locally conformally 
flat manifolds, but this is not necessary in our 
construction. We end with a corollary, whose proof  
is similar to the proof of Theorem \ref{examples}.
\begin{corollary}
\label{examples2}
Let $(M_1,g_1)$ and $(M_2, g_2)$ satisfy
$\int_{M_i} Q_{g_i} dvol_{g_i} \geq 0$, 
and $Y[g_i] > 4 \sqrt{3} \pi$. 
Then the manifold $N = M_1 \# M_2$ admits a metric $\tilde{g}$ 
satisfying (\ref{assump3}). Consequently, $N$ admits 
a metric with $Q = constant.$
\end{corollary}

\bibliography{4dim_JDG_references}
\noindent
\small{\textsc{Department of Mathematics, University of Notre Dame, 
Notre Dame, IN 46556}}\\
{\em{E-mail Address:}} \ {\texttt{mgursky@nd.edu}}\\
\\
\small{\textsc{Department of Mathematics, Massachusetts Institute
of Technology, Cambridge, MA 02139}}\\
{\em{E-mail Address:}} \
{\texttt{jeffv@math.mit.edu}}
\end{document}